\documentclass{amsart}

\usepackage{amsfonts}
\usepackage{amssymb}
\usepackage{amsmath}
\usepackage{graphicx}

\newcommand{\half}{\frac{1}{2}}
\newcommand{\thalf}{\tfrac{1}{2}}

\newcommand{\sump}{\mathop{{\sum}^{+}}}

\def\({\left(}
\def\){\right)}

\numberwithin{equation}{section}

\newtheorem{theorem}{Theorem}[section]
\newtheorem*{theorem*}{Remark}

\newtheorem{lemma}[theorem]{Lemma}

\newtheorem{corollary}[theorem]{Corollary}

\makeatletter
\@namedef{subjclassname@2010}{%
  \textup{2010} Mathematics Subject Classification}
\makeatother

\begin{document}

\title{Nonvanishing of Dirichlet $L$-functions, II}

\author{Rizwanur Khan}
\address{
Mathematics Department, University of Mississippi, University, MS 38677}
\email{rrkhan@olemiss.edu}

\author{Djordje Mili\'cevi\'c}
\address{Bryn Mawr College, Department of Mathematics, 101 North Merion Avenue, Bryn Mawr, PA 19010}
\curraddr{Max-Planck-Institut f\"ur Mathematik, Vivatsgasse 7, D-53111 Bonn, Germany}
\email{dmilicevic@brynmawr.edu}

\author{Hieu T. Ngo}
\address{Vietnam Institute for Advanced Study in Mathematics,
Ta Quang Buu Building, 1 Dai Co Viet Street, Hanoi, Vietnam}
\email{trunghieu.ay@gmail.com}

\thanks{R.K. was supported by the Simons Foundation (award 630985). D.M. was supported by National Science Foundation Grant DMS-1903301.}

\begin{abstract} 
We show that for at least $\frac{5}{13}$ of the primitive Dirichlet characters $\chi$ of large prime modulus, the central value $L(\frac12,\chi)$ does not vanish, improving on the previous best known result of $\frac38$.
\end{abstract}

\maketitle

\section{Introduction}

We revisit the work \cite{khango} of the first and third authors on the nonvanishing of Dirichlet $L$-functions. It was proved there that when $p$ is a large enough prime, for at least three-eighths of the primitive Dirichlet characters $\chi$ of modulus $p$, we have $L(\half, \chi)\neq 0$. For prime moduli (arguably the most interesting case) this was the best result towards the well-known conjecture that $L(\half, \chi)$ should never be zero. In this paper, we introduce two new ingredients which lead to the following improved result.

\begin{theorem}\label{thm:main} Let $\epsilon>0$ be arbitrary. For all primes $p$ large enough in terms of $\epsilon$, for at least $(\frac5{13}-\epsilon)$ of all primitive Dirichlet characters $\chi$ of modulus $p$, we have that $L(\half, \chi)\neq 0$.
\end{theorem}

The nonvanishing problem for Dirichlet $L$-functions has a long history. It was studied by Balasubramanian and Murty \cite{balmur}, who were the first to establish a very small but positive nonvanishing proportion. This was improved significantly by Iwaniec and Sarnak \cite{iwasar} with the proportion $\frac{1}{3}$. Michel and VanderKam \cite{micvan} obtained the same $\frac{1}{3}$ proportion with their symmetric two-piece mollifier described below (as well as nonvanishing results for derivatives of the complete $L$-function). Bui \cite{bui} introduced another type of mollifier to get the proportion 34.11\%. All these results hold for general moduli $q$. For prime moduli $q$, the nonvanishing proportion $\frac{3}{8}$ of \cite{khango} was established by improving for the first time the length of the mollifier from the work of Michel and VanderKam \cite{micvan}. The main novelty in \cite{khango} was to connect the problem to recent advances from the theory of ``trace functions'' (see \cite{fkms} for a wonderful exposition of this new trend of research).

Since the hitherto state of the art paper \cite{khango}, research interest in the nonvanishing of Dirichlet $L$-functions has remained strong. We mention a couple of developments. Bui, Pratt, Robles, and Zaharescu \cite{bprz} were able to increase for the first time the length of the one-piece mollifier introduced by Iwaniec and Sarnak \cite{iwasar}. This has led to some nice applications, but for the nonvanishing problem it does not improve upon \cite{khango}. Pratt \cite{pra} studied the nonvanishing problem on average over the modulus, and obtained a nonvanishing proportion of 50.073\%. We hope that our work will stimulate further research on the nonvanishing of Dirichlet $L$-functions.

{\bf Convention.} Throughout this paper, we adopt the $\epsilon$ convention. That is, $\epsilon$ will denote an arbitrarily small positive constant that may vary from one occurrence to the next. We write $f\ll g$ or $f=O(g)$ to denote $|f|\leqslant Cg$ for some constant $C>0$, which can be made explicit in each instance and may depend on $\epsilon>0$ but is otherwise independent of all parameters including  $p$.

\section{New idea and a proof sketch}\label{background}

We refer the reader to \cite{iwasar,khango,micvan} for background and examples of the method of mollification. In the context of Dirichlet $L$-functions, it involves the evaluation of the first and second mollified moments
\begin{align}
\label{1st} \mathcal{T}_1&=\frac{2}{p} \sump_{\chi \bmod p} L(\thalf, \chi) M_0(\chi),\\
\label{2nd} \mathcal{T}_2&=\frac{2}{p} \sump_{\chi \bmod p} |L(\thalf, \chi) M_0(\chi)|^2,
\end{align}
where the sums are over the even primitive characters (the case of the odd characters being similar). Then an application of the Cauchy--Schwarz inequality shows that $L(\frac12,\chi)\neq 0$ for a proportion of at least $|\mathcal{T}_1|^2/|\mathcal{T}_2|$ of $\chi$. The first idea of mollification is to construct a Dirichlet polynomial $M_0(\chi)$ such that this proportion (which is only $\asymp(\log p)^{-1}$ for the unmollified moments) tends to a positive limit as $p\to\infty$; thus $M_0(\chi)$ is thought of as ``mollifying'' the fluctuations in sizes of central values, and its shape can typically be understood as an instance of a Dirichlet series approximation to $L(s,\chi)^{-1}$ on average.

In \cite{micvan}, Michel and VanderKam used the twisted mollifier given by
\begin{align}
\label{moll} M_0(\chi)= \sum_{m\le M} \frac{y_m \chi(m)}{m^{\half}}+\frac{\overline{\tau}_{\chi}}{p^\half}\sum_{m\le M} \frac{y_m \overline{\chi}(m)}{m^{\half}},
\end{align}
where $M=p^\theta$ for some $\theta>0$, $\tau_\chi$ is the Gauss sum, and
\begin{align*}
y_m=\mu(m)\frac{\log(M/m)}{\log M}.
\end{align*}
It was shown in \cite{micvan} that the above choice for $y_m$ yields a nonvanishing proportion of $\frac{2\theta}{1+2\theta}$. Thus a longer mollifier results in a better nonvanishing proportion. Michel and VanderKam were able to take $\theta<\frac14$, which gave the nonvanishing proportion $\frac13$, while in \cite{khango} Khan and Ngo were able to take $\theta<\frac{3}{10}$, which gave the nonvanishing proportion $\frac38$. 

Our new idea is to replace $M_0(\chi)$ with a more general mollifier of the form
\begin{align}
\label{moll2} M(\chi)= c_1\sum_{m\le M R} \frac{y_m \chi(m)}{m^{\half}}+c_2\frac{\overline{\tau}_{\chi}}{p^\half}\sum_{m\le M } \frac{y_m \overline{\chi}(m)}{m^{\half}},
\end{align}
where $M=p^\theta$ as before, $R=p^\alpha$ for some $\alpha>0$, and $c_1,c_2>0$. (Equivalently, we could take one sum of length $p^{\theta_1}$ and the other of length $p^{\theta_2}$.) This is novel in two related but distinct ways. First off, we are using an unbalanced mollifier comprising of sums of unequal lengths, in a way reminiscent of using an unbalanced functional equation to represent an $L$-function. Second, the two pieces of the mollifier are assigned different weights, which (keeping in mind the paradigm of $L(s,\chi)^{-1}$ and the duality between $M(\chi)$ and $M(\chi^{-1})$) can also be thought of as normalizing weights accounting for the length. We will comment further on our choice of the mollifier later in the introduction as well as in Section \ref{sect:optimization}.

Our goal is to evaluate the first and second mollified moments
\begin{align}
\label{1stnew} \mathcal{S}_1=&\frac{2}{p} \sump_{\chi \bmod p} L(\thalf, \chi) M(\chi),\\
\label{2ndnew} \mathcal{S}_2=&\frac{2}{p} \sump_{\chi \bmod p} |L(\thalf, \chi) M(\chi)|^2.
\end{align}

The mollified second moments $\mathcal{T}_2$ and $\mathcal{S}_2$ are harder to treat than the mollified first moments $\mathcal{T}_1$ and $\mathcal{S}_1$. Let 
$$S(x,y;p)=\sum_{\substack{u \bmod p\\ u\overline{u} \equiv 1 \bmod p}} e\Big(\frac{xu+y\overline{u}}{p}\Big)$$
denote the Kloosterman sum. In \cite{micvan}, when \eqref{moll} is substituted into \eqref{2nd} and the square expanded, the cross terms involving the Gauss sum are the hardest to treat. 
After some tranformations involving Poisson summation, these cross terms lead to the problem of proving an estimate of the shape
\begin{align*}
\frac{1}{pM^2} \sum_{n,k,m_1,m_2\sim M} y_{m_1} y_{m_2}S(nk\overline{m}_1,\overline{m}_2;p)\ll p^{-\epsilon}.
\end{align*}
When $M=p^\frac{1}{4}$, the problem is precisely to break the ``trivial'' estimate given by inserting Weil's bound for Kloosterman sums. This seems like a natural barrier. To do better and take $M>p^\frac14$ in at least some of the variables, one must be able to detect cancellation between the Kloosterman sums. The strategy in \cite{khango} was to glue together some of the variables, writing $h=nk\overline{m}_1$, and then to apply H\"{o}lder's inequality, getting 
\begin{equation}
\label{holder}
\begin{aligned}
&\frac{1}{pM^2} \sum_{n,k,m_1,m_2\sim M} y_{m_1} y_{m_2} S(nk\overline{m}_1,\overline{m}_2;p) \\ 
&\qquad \ll \frac{p^\epsilon}{pM^2} \Big( \sum_{h\bmod p} \nu(h)^\frac43 \Big)^{\frac{3}{4}} \Big( \sum_{h\bmod p}  \Big|\sum_{m_2\sim M} y_{m_2} S(h,\overline{m}_2;p)\Big|^4 \Big)^{\frac{1}{4}} ,
\end{aligned}
\end{equation}
where $\nu(h)$ denotes the number of ways of writing $h$ as $nk\overline{m}_1$ mod $p$. If $M^3 \le p^{1-\epsilon}$ then 
\begin{equation}
\label{mult}
\begin{aligned}
\sum_{h\bmod p} \nu(h)^\frac43\ll \sum_{h\bmod p} \nu(h)^2
&\ll \sum_{nk\overline{m}_1\equiv n'k'\overline{m}_1'\bmod p} 1\\
&\ll \sum_{nk{m}_1'\equiv n'k'{m}_1\bmod p} 1 = \sum_{nk{m}_1'=n'k'{m}_1} 1 \ll M^{3+\epsilon}.
\end{aligned}
\end{equation}
Thus while the first $h$-sum in \eqref{holder} counts about $M^3$ elements, the second $h$-sum has been extended to a complete sum mod $p$. This is a wasteful step, but not too bad because $M^3$ is at least $p^\frac{3}{4}$ in the ranges of interest. In this way, we get the upper bound
\begin{equation}
\label{sum-of-sums-of-products}
\frac{p^\epsilon}{pM^2} (M^3)^\frac34 \Big( \sum_{m_1,m_2,m_3,m_4\sim M} \Big| \sum_{h\bmod p} S(h, \overline{m}_1;p)S(h, \overline{m}_{2};p)S(h ,\overline{m}_{3};p)S(h, \overline{m}_{4};p)\Big| \Big)^{\frac{1}{4}}.
\end{equation}
For the innermost $h$-sum, we get from \cite[Proposition 3.2]{fgkm} that if at least one number in the tuple $(m_1,m_2,m_3,m_4)$ is distinct from the others, then $\sum_{h\bmod p}\ll p^\frac52 $, which saves a factor $p^\frac12 $ over Weil's bound. The other case where no numbers in the tuple are distinct from the others forms a smaller set for which we can just apply Weil's bound. The limitation of this method is $M=p^\frac{3}{10}$. It was not explicitly shown in \cite{khango}, but other choices of H\"{o}lder exponents do not yield good results.

With our new mollifier $M(\chi)$, the crucial estimate to show will be roughly
\begin{equation}
\label{roughly-to-estimate}
\frac{1}{pM^2R} \sum_{\substack{n,k\sim M\sqrt{R} \\ m_2 \sim M\\ m_1\sim MR}} y_{m_1} y_{m_2}S(nk\overline{m}_1,\overline{m}_2;p)\ll p^{-\epsilon}.
\end{equation}
The point is that when we glue together $h=nk\overline{m}_1$, the unbalanced mollifier allows us to use the flexibility to increase the length of $m_1$, and therefore have $h$ cover more elements mod $p$. This way, when we extend the sum over $h$ to a complete sum mod $p$, this is less wasteful (or not at all, up to $p^{\epsilon}$ factors, if the method can be pushed sufficiently far). However, while the additional freedom certainly does not hurt, it is not clear a priori whether or not it improves the result because as we take $R$ larger, we may need to take $M$ smaller. That is, when we make one of the sums comprising the mollifier longer, the other sum may need to be shorter in order to maintain control of the error terms. 

Specifically, using the unbalanced mollifier and starting from \eqref{roughly-to-estimate} and proceeding as above, we find in Burgess' trick \eqref{mult} that if $M^3R^2\leq p^{1-\epsilon}$, the sum on the left-hand side is $\ll (M^3R^2)^{1+\epsilon}$, and thus estimating the complete exponential sums at the point of \eqref{sum-of-sums-of-products} we find that the left-hand side of \eqref{roughly-to-estimate} is
\begin{equation}
\label{after-estimation}
\ll\frac{p^{\epsilon}}{pM^2R}(M^3R^2)^{\frac34}\big(M^2p^3+M^4p^{\frac52}\big)^{\frac14}=\frac1{p^{\frac14-\epsilon}}M^{\frac34}R^{\frac12}+\frac1{p^{\frac38-\epsilon}}M^{\frac54}R^{\frac12}.
\end{equation}

To see how well our unbalanced mollifier \eqref{moll2} works altogether, we fast forward to the end game; the main terms are treated as in \cite{khango}, bookkeeping for weights and mollifier lengths, and we find that
\[ \mathcal{S}_1\sim c_1+c_2,\quad \mathcal{S}_2\sim (c_1+c_2)^2+c_1^2/(\theta+\alpha)+c_2^2/\theta \]
as $p\to\infty$, subject to the conditions $3\theta+2\alpha<1$ and $10\theta+4\alpha<3$ arising from the estimation \eqref{after-estimation}. Optimizing the ratio $|\mathcal{S}_1|^2/|\mathcal{S}_2|$ leads to the choices $M=p^{\frac14-\epsilon}$, $R=p^{\frac18-\epsilon}$, $c_1=\frac35$, $c_2=\frac25$. We note that our choices $\theta$ and $\alpha$ are essentially as large as possible; the normalizing weights $c_1$ and $c_2$ can be optimized subsequently to $c_1:c_2=\log(MR):\log(M)$. In particular, the passage to the complete sum over $h$ modulo $p$ in \eqref{holder} is essentially saturated, so in principle costless up to $p^{\epsilon}$ factors.

One way to understand parameters $c_1$, $c_2$ is that, in an unbalanced amplifier, the longer part does a better job at mollification, so we assign it a higher weight. Alternatively, optimizing a general unbalanced mollifier using coefficients $(z_m)$ and $(z'_m)$ in the direct and dual sums in \eqref{moll2} leads to maximizing a quadratic form in $(z_m)$, $(z'_m)$ subject to a linear constraint, with the same optimal choice $z_m/c_1=z'_m/c_2=y_m$. Finally we remark that, in all approaches including \cite{iwasar,micvan,khango} and ours, the proportion of nonvanishing obtained equals $\vartheta/(1+\vartheta)$, where $p^{\vartheta}$ is the \textit{combined} length of the mollifier ($MR\cdot M$ in our case). From this perspective, the improved treatment of the cross terms allows for increase from $\vartheta<\frac12$ of \cite{micvan} to the present $\vartheta<\frac58$, with the full increase of the combined length allocated to a single piece of the optimal mollifier.

\section{Mollified moments}

In this section we describe how to apply the mollification process of \cite{micvan} and \cite{khango} with the new mollifier \eqref{moll2}. Both these works use the mollifier \eqref{moll2} with $R=1$ (equivalently, $\alpha=0$), so that the sums have equal length $M=p^\theta$. The main terms of the mollified first and second moments $\mathcal{T}_1$ and $\mathcal{T}_2$ yield a nonvanishing proportion of $\frac{2\theta}{1+2\theta}$. We first examine the mollification to find the corresponding proportion of nonvanishing when a mollifier with unequal length sums is used instead.

\subsection{Mollified first moment}

We first recall the approximate functional equation \cite[Theorem 5.3]{iwakow}
\begin{align}
\label{afe2} L(\thalf, \chi) = \sum_{n\ge 1} \frac{\chi(n)}{n^\half}W\Big(\frac{n}{p^\half}\Big) + \frac{\tau_\chi}{p^\half} \sum_{n\ge 1} \frac{\overline{\chi}(n)}{n^\half}W\Big(\frac{n}{p^\half}\Big),
\end{align}
where
\begin{align*}
W(x)=\frac{1}{2\pi i} \int_{(2)} \frac{\Gamma(\frac{s}{2}+\frac{1}{4})}{\Gamma( \frac{1}{4})}(\pi^\half x)^{-s} \frac{ds}{s}.
\end{align*}
Inserting \eqref{moll2} and \eqref{afe2} into \eqref{1stnew} and expanding, we see that the first mollified moment equals
\begin{align}
\label{s1expand} \mathcal{S}_1= \frac{2}{p} \sump_{\chi \bmod p}\Bigg( &\sum_{\substack{n\ge 1\\m\le M R}} \frac{c_1 y_m \chi(nm)}{(nm)^{\half}} W\Big(\frac{n}{p^\half}\Big) +\sum_{\substack{n\ge 1\\m\le M} } \frac{c_2 y_m \overline{\chi}(nm)}{(nm)^{\half}} W\Big(\frac{n}{p^\half}\Big)  \\
\nonumber +\frac{\overline{\tau}_\chi}{p^\half}&\sum_{\substack{n\ge 1\\m\le M }} \frac{c_2 y_m \chi(n)\overline{\chi}(m)}{(nm)^{\half}} W\Big(\frac{n}{p^\half}\Big) +\frac{\tau_\chi}{p^\half}\sum_{\substack{n\ge 1\\m\le M R}} \frac{c_1 y_m \chi(m)\overline{\chi}(n)}{(nm)^{\half}} W\Big(\frac{n}{p^\half}\Big) \Bigg).
\end{align}
A similar expression of course holds for $\mathcal{T}_1$, the asymptotic evaluation of which is given in \cite[equation (9)]{micvan}. This proceeds by treating each of the sums as in \eqref{s1expand} separately; therefore we can directly import the results of \cite{micvan} into \eqref{s1expand}, with some bookkeeping to account for $c_1$ and $c_2$, lengths of the $m$-sums, and the fact that Michel and VanderKam work with the completed Dirichlet $L$-functions and use slightly different normalizations, cf. \eqref{afe2}, \eqref{s1expand}, and \cite[p.131 and equation (2)]{micvan}. The main term arises from the $n=m=1$ terms of the first line of \eqref{s1expand}, which contribute $c_1+c_2+O(p^{-\frac14+\epsilon})$. The error term in this evaluation is $O(p^{-\epsilon})$ as long as the sums comprising the mollifier have lengths less than $p^{\half-\epsilon}$, as stated right above \cite[equation (9)]{micvan}, keeping in mind the notation $\hat{q}=(q/\pi)^\half$ of \cite{micvan}, so that $\Delta$ translates to $2\theta$ for a mollifier of length $p^{\theta}$. The condition $MR\leq p^{\frac12-\epsilon}$ will be ensured by the first two conditions of \eqref{eq:optimization-constraint}, and so
\begin{align}
\label{s1final} \mathcal{S}_1=c_1+c_2+O(p^{-\epsilon}).
\end{align}

\subsection{Mollified second moment}
For the second mollified moment, inserting \eqref{moll2} into \eqref{2ndnew} and expanding, we see that $\mathcal{S}_2$ given in \eqref{2ndnew} equals
\begin{align}
\label{2nd1} &2c_1c_2\frac{2}{p}\sump_{\chi\bmod p} |L(\thalf,\chi)|^2 \frac{{\tau}_{\chi}}{p^\half} \sum_{\substack{m_1 \le M R\\ m_2\le M}}   \frac{y_{m_1}y_{m_2} \chi(m_1)\chi(m_2)}{(m_1m_2)^{\half}}\\
\label{2nd2}&\qquad+c_1^2\frac{2}{p} \sump_{\chi\bmod p} |L(\thalf,\chi)|^2 \Big|\sum_{m\le M R} \frac{y_m \chi(m)}{m^{\half}}\Big|^2 + c_2^2\frac{2}{p} \sump_{\chi\bmod p} |L(\thalf,\chi)|^2 \Big|\sum_{m\le M } \frac{y_m \chi(m)}{m^{\half}}\Big|^2.
\end{align}
The standard first step in the evaluation of \eqref{2nd1} and \eqref{2nd2} comprises of firstly proving the approximate functional equation (see \cite[equation (2--2)]{khango})
\begin{align}
\label{afe} |L(\thalf, \chi)|^2= 2\sum_{n_1,n_2\ge 1} \frac{\chi(n_1)\overline{\chi}(n_2)}{(n_1n_2)^\half}V\Big(\frac{n_1n_2}{p}\Big)
\end{align}
where
\begin{align*}
V(x)=\frac{1}{2\pi i} \int_{(2)} \frac{\Gamma(\frac{s}{2}+\frac{1}{4})^2}{\Gamma( \frac{1}{4})^2}(\pi x)^{-s} \frac{ds}{s},
\end{align*}
secondly inserting \eqref{afe} into \eqref{2nd1} and \eqref{2nd2}, and thirdly invoking for $(n_1n_2,p)=1$ the approximate identities (see \cite[equation (17)]{micvan} or \cite[equation (3.4)]{iwasar}),
\begin{align*}
&\frac{2}{p} \sump_{\chi\bmod p}  \chi(n_1)\overline{\chi}(n_2) = 
\begin{cases}
1 + O(p^{-1}) &\text{ if } n_1\equiv \pm n_2 \,\, ({\rm mod}\,\, p)\\
O(p^{-1}) &\text{ otherwise,}
\end{cases} \\
&\frac{1}{p} \sump_{\chi\bmod p} \tau_\chi \chi(n_1) =  \text{Re}\Big(e\Big(\frac{\overline{n_1}}{p}\Big)\Big)+O(p^{-1}), \ \ \ \  
\end{align*}
where as usual $e(x)=e^{2\pi i x}$. The output will be representations of \eqref{2nd1} and \eqref{2nd2} as quadruple sums which we then separate into main terms and error terms. The smooth function $V$ has the effect of imposing the condition $n_1n_2\leq p^{1+\epsilon}$ on the variables $n_1$ and $n_2$, because on moving the line of integration we infer that $V(x)\ll_c x^{-c}$ for any $c>0$. 

For \eqref{2nd2}, we look to \cite[section 5]{micvan}, where sums of the type 
\begin{align}
\label{type} \frac{2}{p} \sump_{\chi\bmod p}  |L(\thalf,\chi)|^2 \Big|\sum_{m\le p^{\theta}} \frac{y_m \chi(m)}{m^{\half}}\Big|^2
\end{align}
are evaluated starting from the standard first step just described. Specifically, referring to \cite[equation (16)]{micvan} for the main term and \cite[page 136]{micvan} for the error term, we get that \eqref{type} equals
$\frac1{\theta}+1+ O(p^{-\epsilon})$ for $\theta<\half$. This can be seen by taking $k=0$ and $P_0(t)=t$ in \cite[equation (16)]{micvan}, removing the factor $\Gamma(\frac14)^2 \hat{q} \phi^+(q)$ therein (which, as in the first moment, is due to the use of completed $L$-functions and slightly different normalizations in \cite{micvan}; cf.~\eqref{afe}, \eqref{type}, and \cite[p.131 and equation (3)]{micvan}), and keeping in mind that $\Delta=2\theta$. Using this evaluation of \eqref{type}, we get that \eqref{2nd2} equals
\[ c_1^2\Big(\frac{1}{\theta+\alpha}+1\Big) + c_2^2\Big(\frac{1}{\theta}+1\Big) + O(p^{-\epsilon})  \]
for $\theta+\alpha<\half.$ 

We now consider \eqref{2nd1}. The above standard first step implies that \eqref{2nd1} equals
\begin{align}
\label{2moll} \frac{4c_1c_2}{p^\half} {\rm Re} 
\sum_{\substack{n_1,n_2\ge 1\\ m_1\le MR \\ m_2\le M \\(n_1 n_2 m_1m_2,p)=1}} 
\frac{y_{m_1}y_{m_2}}{(n_1n_2m_1m_2)^\half} V\Big(\frac{n_1n_2}{p}\Big) e\Big(\frac{ n_2 \, \overline{n_1 m_1 m_2}}{p}\Big) + O\Big(\frac{MR^{\frac{1}{2}}}{p^{1-\epsilon}}\Big).
\end{align}
It is shown in \cite[section 6.1]{micvan} that \eqref{2moll} yields a constant main term which is contained in the contribution of the terms with $m_1m_2n_1=1$. Specifically, it follows from \cite[equation (23)]{micvan}, after again removing the normalization factor $\Gamma(\frac14)^2 \hat{q} \phi^+(q)$ (of the same origin as above), that the contribution of these terms to
\eqref{2moll} is
\[ 2c_1c_2+O(p^{-\epsilon}), \]
where it is easy to see that the method indeed gives power savings for $k=0$.

Next, we consider the terms with $m_1m_2n_1>1$ in dyadic intervals.
Let
\begin{equation}\label{bilinear-term-predefinition}
\mathcal{B}(M_1,M_2,N_1,N_2)=
	\sum_{\substack{N_1\le n_1 \le 2N_1 \\ N_2\le n_2\le 2N_2\\
		M_1\le m_1< 2M_1 \\ M_2\le m_2< 2M_2}} 
	\frac{y_{m_1} y_{m_2} }{(pM_1M_2N_1N_2)^\half} V\Big(\frac{n_1n_2}{p}\Big) 
	f_1\Big(\frac{n_1}{N_1}\Big) f_2\Big(\frac{n_2}{N_2}\Big) 
	e\Big(\frac{n_2\overline{n_1 m_1 m_2}}{p}\Big)  ,
\end{equation}
for any integers $M_1,M_2\ge 1$ and any $N_1, N_2\ge \frac12$ satisfying
\begin{align}
\label{ranges} N_1N_2\le p^{1+\epsilon}, \ \ \ M_1 \le \frac{MR+1}{2}, \ \ \ M_2 \le \frac{M+1}{2}, \ \ \ M_1M_2N_1\ge 1,
\end{align}
arbitrary coefficients $y_{m_1}, y_{m_2}$ and any fixed smooth functions $f_1,f_2$ compactly supported on the interval $(\frac54,\frac74)$ say, all with absolute values bounded by $p^\epsilon$.
We remark that the definition of $\mathcal{B}(M_1,M_2,N_1,N_2)$ is slightly different from \cite{khango}; still we choose the same notation because of the similarity. Note that the last condition of \eqref{ranges} restricts the sum to $m_1m_2n_1>1$. This is because if $M_1M_2=1$ then this condition implies that $N_1\ge 1$, in which case the support of the test functions implies that $f_1(\frac{n_1}{N_1})$ vanishes unless $n_1>1$.

On putting 
$f(n_1,n_2)= V(\frac{n_1n_2}{p}) f_1(\frac{n_1}{N_1}) f_2(\frac{n_2}{N_2})$, we can rewrite 
\begin{equation}\label{bilinear-term-definition}
\mathcal{B}(M_1,M_2,N_1,N_2)= \frac{1}{(pM_1M_2N_1N_2)^\half} \sum_{\substack{N_1\le n_1 \le 2N_1 \\ N_2\le n_2\le 2N_2\\M_1\le m_1< 2M_1 \\ M_2\le m_2< 2M_2}} y_{m_1} y_{m_2} f(n_1,n_2) e\Big(\frac{n_2\overline{n_1 m_1 m_2}}{p}\Big)  ,
\end{equation}
with coefficients $y_{m_1}, y_{m_2}$ and the smooth function $f$ all having absolute values bounded by $p^\epsilon$.
In Section \ref{sect:error}, specifically Corollary \ref{cor:bilinear-bound}, we shall establish the bound
\begin{equation}\label{bilinear-term-bound}
\mathcal{B}(M_1,M_2,N_1,N_2) \ll p^{-\epsilon}
\end{equation}
for $M=p^\theta$, $R=p^\alpha$ with 
\begin{equation}\label{eq:optimization-constraint}
0<\theta<\frac{1}{2}-\epsilon, \ \ \ \ 0<\theta+\alpha<\frac{1}{2}-\epsilon, \ \ \ \ 
3\theta+2\alpha-1<0, \ \ \ \ 
10\theta+4\alpha-3<0.
\end{equation}
It will thus follow that, under the assumption \eqref{eq:optimization-constraint}, the terms with $m_1m_2n_1>1$ in \eqref{2moll} contribute $O(p^{-\epsilon})$. These conditions also ensure that each component of our two-piece mollifier has length less than $p^{\half-\epsilon}$. 

In summary, combining our evaluations of \eqref{2nd1} and \eqref{2nd2}, we find that 
\begin{align}
\label{s2final} \mathcal{S}_2 &=c_1^2\Big(\frac{1}{\theta+\alpha}+1\Big) + c_2^2\Big(\frac{1}{\theta}+1\Big) + 2c_1c_2+O(p^{-\epsilon}) \\
\nonumber &= \frac{c_1^2}{\theta+\alpha}+\frac{c_2^2}{\theta}+(c_1+c_2)^2+O(p^{-\epsilon}),
\end{align}
assuming \eqref{eq:optimization-constraint}.

\subsection{Nonvanishing proportion} We are in a position to derive the proportion of nonvanishing as a function of the lengths of the mollifier components and weights in \eqref{moll2}.

\begin{lemma}\label{lem:newprop}
Using the mollifier \eqref{moll2} and assuming the condition \eqref{eq:optimization-constraint}, the main terms of the mollified first and second moments given in \eqref{1st} and \eqref{2nd} yield the nonvanishing proportion
\begin{align}
\label{prop} \Big(\frac{(c_1/(c_1+c_2))^2}{\theta+\alpha}+\frac{(c_2/(c_1+c_2))^2}{\theta}+1\Big)^{-1}.
\end{align}
\end{lemma}

\proof
As is standard, and mentioned in the introduction, the nonvaninishing proportion is at least $|\mathcal{S}_1|^2/|\mathcal{S}_2|$. We now insert \eqref{s1final} and \eqref{s2final} to complete the proof. 
\endproof

\section{Error term}\label{sect:error}

In this section we set out to prove the estimate \eqref{bilinear-term-bound} for the sum $\mathcal{B}(M_1,M_2,N_1,N_2)$ given by \eqref{bilinear-term-definition} under the condition \eqref{ranges}. Now there are two natural ways to proceed. 
In \eqref{bilinear-term-definition}, on applying Poisson summation in $n_2$ after first separating into residue classes modulo $p$, we get the following estimate.
\begin{lemma}\label{pois1} For $MR^2<p^{1-\epsilon}$, we have
\begin{align*}
\mathcal{B}(M_1,M_2,N_1,N_2) \ll p^\epsilon \Big(\frac{M_1M_2 N_1}{pN_2}\Big)^{1/2} + p^{-\epsilon} \ll p^\epsilon \Big(\frac{M^2 R N_1}{pN_2}\Big)^{1/2} + p^{-\epsilon}.
\end{align*}
\end{lemma}
\proof 
This is given by \cite[equation (27)]{micvan}, or equivalently \cite[equation (2-6)]{khango}.
We note that this is the bound that would not cover the summands with $m_1m_2n_1=1$ (cf.~\cite[p.146]{micvan}), but that possibility is excluded in the definition of $\mathcal{B}(M_1,M_2,N_1,N_2)$.
\endproof

In \eqref{bilinear-term-definition}, if we instead separate $n_1$ into residue classes modulo $p$ and apply Poisson summation, denoting the dual variable by $k$, then we get Kloosterman sums as follows.

\begin{lemma}\label{lem:bilinear} Suppose $MR^2<p^{1-\epsilon}$. For some function $\hat{f}$ with $\| \hat{f} \|_\infty \ll p^\epsilon$, we have
\begin{align}
\label{newb} \mathcal{B}(M_1,M_2,N_1,N_2)= 
	\frac{1}{(pM_1M_2N_1N_2)^\half} \frac{N_1}{p} 
	\sum_{\substack{ 1\le |k| \le p^{1+\epsilon}/N_1\\ N_2\le n_2\le 2 N_2 \\
	M_1\le m_1< 2M_1 \\ M_2\le m_2< 2M_2}} 
	y_{m_1} y_{m_2} \hat{f}(n_2,k) S(kn_2,\overline{m_1m_2},p) + O(p^{-\epsilon}).
\end{align}
\end{lemma}
\proof
See \cite[equation (3-5)]{khango}. The contribution of $k=0$ is shown on \cite[page 8]{khango} to be $O(p^{-\epsilon})$ when $MR^2<p^{1-\epsilon}$.
\endproof

We can estimate \eqref{newb} as follows.

\begin{lemma}\label{lem:trilinear}
Suppose $\frac{p}{N_1}N_2 M_1<p^{1-\epsilon}$. We have
$$
\sum_{\substack{1\leq |k| \leq p^{1+\epsilon}/N_1 \\ N_2 \leq n_2 \leq 2N_2 \\ 
	M_1\le m_1< 2M_1 \\ M_2\le m_2< 2M_2}}
y_{m_1}y_{m_2}\hat{f}(n_2,k)S(kn_2\overline{m}_1, \overline{m}_2;p)
\ll 
p^\epsilon \Big( \frac{pN_2M_1}{N_1} \Big)^{\frac{3}{4}} (M_2p^{\frac58} + M_2^{\frac12}p^{\frac34}).
$$
\end{lemma}
\proof
This is an immediate consequence of \cite[Lemma 3.2]{khango}, but we review the proof. We glue together $h = k n_2 \overline{m}_1$ and apply H\"{o}lder's inequality as described in Section \ref{background}.
Provided $\frac{p}{N_1}N_2 M_1<p^{1-\epsilon}$, we have by the argument in \eqref{mult} that 
\begin{multline*}
\sum_{\substack{1\leq |k| \leq p^{1+\epsilon}/N_1 \\ N_2 \leq n_2 \leq 2N_2 \\ 
	M_1\le m_1< 2M_1 \\ M_2\le m_2< 2M_2}}
y_{m_1}y_{m_2}\hat{f}(n_2,k)S(kn_2\overline{m}_1, \overline{m}_2;p)\\
\ll  p^\epsilon \Big( \frac{pN_2M_1 }{N_1} \Big)^\frac34
  \Big( \sum_{m_1,m_2,m_3,m_4<2M_2} \Big| \sum_{h\bmod p}S(h,\overline{m}_1;p)S(h,\overline{m}_{2};p)S(h,\overline{m}_{3};p)S(h,\overline{m}_{4};p)\Big| \Big)^{\frac{1}{4}}.
\end{multline*}
 The number of tuples $(m_1,m_2,m_3,m_4)$ where no entry is distinct from the others is $O(M_2^2)$. To these tuples we apply Weil's bound for Kloosterman sums. For the rest, we use \cite[Proposition 3.2]{fgkm} to get that $\sum_{h \bmod p}\ll p^\frac52$, and deduce
the lemma.
\endproof

Putting Lemma \ref{lem:bilinear} and Lemma \ref{lem:trilinear} together, keeping in mind $M_1\le MR$ and $M_2\le M$, we get

\begin{lemma} \label{bbound} Suppose $MR^2<p^{1-\epsilon}$ and $\frac{p}{N_1}N_2 M_1<p^{1-\epsilon}$. We have
\begin{align*}
\mathcal{B}(M_1,M_2,N_1,N_2) 
\ll \Big(\frac{p^\epsilon N_2 MR}{N_1}\Big)^\frac14 +  \Big(\frac{p^\epsilon N_2^2M^6R^2}{N_1^2p}\Big)^\frac18 + O(p^{-\epsilon}).
\end{align*}
\end{lemma}

Finally, we are ready to prove

\begin{corollary}\label{cor:bilinear-bound}
Under the assumption \eqref{eq:optimization-constraint}, we have that $\mathcal{B}(M_1,M_2,N_1,N_2) \ll  p^{-\epsilon}$.
\end{corollary}

\proof
We consider two cases. If $\frac{M^2 R N_1}{pN_2}<p^{-\epsilon}$, then we are done by Lemma \ref{pois1}. Note that the condition of Lemma \ref{pois1} is satisfied by \eqref{eq:optimization-constraint}. Therefore assume
\begin{align}
\label{assume} \frac{M^2 R N_1}{pN_2}\geq p^{-\epsilon}.
\end{align}
We can combine this assumption with \eqref{eq:optimization-constraint} to see that the conditions of Lemma \ref{bbound} are satisfied. Also, using \eqref{assume}, we see that the estimate of Lemma \ref{bbound} implies that
\begin{align*}
\mathcal{B}(M_1,M_2,N_1,N_2)\ll \Big(\frac{p^\epsilon M^3 R^2}{p}\Big)^\frac14 +  \Big(\frac{p^\epsilon M^{10} R^4}{p^3}\Big)^\frac18 + O(p^{-\epsilon}).
\end{align*}
It remains to observe that the assumption \eqref{eq:optimization-constraint} is precisely what makes $\mathcal{B}(M_1,M_2,N_1,N_2) \ll  p^{-\epsilon}$. The corollary is proved.
\endproof

\section{Optimization}\label{sect:optimization}

We are in a position to prove our main theorem.

\proof[Proof of Theorem \ref{thm:main}]
In view of Lemma \ref{lem:newprop}, the task now is to minimize the ratio
\begin{equation}
\label{to-minimize}
\Big(\frac{c_1^2}{\theta+\alpha}+\frac{c_2^2}{\theta}\Big):(c_1+c_2)^2
\end{equation}
subject to the conditions $c_1,c_2>0$ as well as
\begin{align}
\label{cond0} &0<\theta<\frac{1}{2}-\epsilon, \ \ \ \ 0<\theta+\alpha<\frac{1}{2}-\epsilon,\\
\label{cond1} &3\theta+2\alpha-1<0,\\
\label{cond2} &10\theta+4\alpha-3<0.
\end{align}
These are the conditions of Corollary \ref{cor:bilinear-bound}. Recall that the condition \eqref{cond0} guarantees that the lengths of the mollifier components do not exceed $p^{\frac{1}{2}-\epsilon}$, whereas the conditions \eqref{cond1} and \eqref{cond2} ensure that $\mathcal{B}(M_1,M_2,N_1,N_2)\ll p^{-\epsilon}$ by Lemma \ref{bbound}.

Take for a moment arbitrary $\theta,\alpha>0$ satisfying the conditions \eqref{cond0}--\eqref{cond2}. Then it is immediate from the Cauchy--Schwarz inequality, applied to the sum $\sqrt{\theta+\alpha}\frac{c_1}{\sqrt{\theta+\alpha}}+\sqrt{\theta}\frac{c_2}{\sqrt{\theta}}$, that the ratio \eqref{to-minimize} is minimized when $c_1:c_2=(\theta+\alpha):\theta$, in which case we find from Lemma~\ref{lem:newprop} that we obtain the proportion of nonvanishing
\[ \Big(1+\frac1{2\theta+\alpha}\Big)^{-1}. \]
Note that, as alluded to in the introduction, this expression has a natural intrinsic meaning, namely $p^{2\theta+\alpha}=MR\cdot M$ is the combined length of the mollifier \eqref{moll2}.

Thus the question is to maximize this length subject to the conditions \eqref{cond0}--\eqref{cond2}. Under those conditions, we have that
\[ 2\theta+\alpha=\frac14(3\theta+2\alpha)+\frac18(10\theta+4\alpha)<\frac58, \]
and the choice $\theta=\frac14-\epsilon$, $\alpha=\frac18-\epsilon$ is essentially optimal, verifying \eqref{cond0}--\eqref{cond2} and giving the announced proportion of nonvanishing
\[ \Big(1+\frac{1}{5/8}\Big)^{-1}-\epsilon=\frac5{13}-\epsilon. \qedhere \]
\endproof

Note that the sums comprising our mollifier can have lengths up to $MR=p^{\frac38-\epsilon}$ and $M=p^{\frac14-\epsilon}$. In comparison with \cite{khango}, which essentially corresponds to \eqref{cond0}--\eqref{cond2} with $\alpha=0$ and had both components of the mollifier of length $p^{\frac{3}{10}-\epsilon}$, the direct (longer) sum in \eqref{moll2} is now quite a bit longer but the dual (shorter) sum a bit shorter, with the increased combined length responsible for the improved proportion of nonvanishing.

It is also instructive to recall the earlier work of Michel--VanderKam~\cite{micvan}, which added the flexibility of the two-piece mollifier, with both components of length $p^{\frac14-\epsilon}$, the maximum allowed from the Weil bound treatment of the cross terms. The improved treatment of cross terms from \cite{khango} coupled with the unbalanced mollifier of the present paper gives rise to the substantially relaxed conditions \eqref{cond1}--\eqref{cond2}, which are now both saturated, and where the entire gain from the increased combined length (from $p^{\frac12-\epsilon}$ up to $p^{\frac58-\epsilon}$) is optimally attributed to one of the components of the mollifier \eqref{moll2}.

\bibliographystyle{amsplain}

\bibliography{nonvanishing-part2j}

\end{document}